\documentstyle[10pt]{amsart}

\newtheorem{prop}{Proposition}

\newtheorem{lemma}{Lemma}
\newtheorem{cor}{Corollary}

\title{Twistorial construction of minimal hypersurfaces}
\author{Johann Davidov}

\address{Institute of Mathematics and Informatics \\
Bulgarian Academy of Sciences\\ Acad. G.Bonchev Str. Bl.8\\ 1113 Sofia\\ Bulgaria and \\ "L.
Karavelov" Civil Engineering Higher School\\ 1373 Sofia.}\email{jtd@@math.bas.bg}

\begin{document}

\begin{abstract}
Every almost Hermitian structure $(g,J)$ on a four-manifold $M$ determines a hypersurface
$\Sigma_J$ in the (positive) twistor space of $(M,g)$ consisting of the complex structures
anti-commuting with $J$. In this note we find the conditions under which $\Sigma_J$ is minimal with
respect to a natural Riemannian metric on the twistor space in the cases when $J$ is integrable or
symplectic. Several examples illustrating the obtained results are also discussed.

\vspace{0,1cm} \noindent {\it Keywords:} Twistor spaces; minimal hypersurfaces.

\vspace{0,1cm} \noindent {\it Mathematics Subject Classification 2010}. Primary: 53C28; Secondary:
53A10, 49Q05.

\end{abstract}

\thispagestyle{empty}

\maketitle
\vspace{0.5cm}

\section{Introduction}

The twistor space ${\cal Z}$ of a Riemannian manifold $(M,g)$ is the bundle on $M$ parametrizing
the complex structures on the tangent spaces of $M$ compatible with the metric $g$. Thus the almost
Hermitian structures on $(M,g)$ are sections of ${\cal Z}$. Given such a structure $J$, we can
consider the hypersurface $\Sigma_J$  of points of ${\cal Z}$ representing complex structures
anti-commuting with $J$.  The twistor space admits a $1$-parameter family $h_t$ of Riemannian
metrics, the so-called canonical variation of $g$. Then it is natural to relate geometric
properties of the hypersurface $\Sigma_J$ in the Riemannian manifold $({\cal Z},h_t)$ to properties
of the almost Hermitian structure $(g,J)$. In this note we address the problem of when $\Sigma_J$
is a minimal hypersurface in the twistor space of a manifold of dimension four. In this dimension,
there are three basic classes in the Gray-Hervella classification - those of Hermitian, almost
K\"ahler (symplectic) and K\"ahler manifolds. If $(g,J)$ is K\"ahler, $\Sigma_J$  is a totally
geodesic submanifold, as one can expect. In the case of an Hermitian manifold, we express the
condition for minimality of $\Sigma_J$ in terms of the Lee form of $(M,g,J)$, while for an almost
K\"ahler manifold we show that $\Sigma_J$ is minimal if and only if the $\star$-Ricci tensor of
$(M,g,J)$ is symmetric. Several example illustrating these results are discussed in the last
section of the paper.

\section{Preliminaries}

Let $(M,g)$ be an oriented  Riemannian manifold of dimension four.
The metric $g$ induces a metric on the bundle of two-vectors
$\pi:\Lambda^2TM\to M$ by the formula
$$
g(v_1\wedge v_2,v_3\wedge v_4)=\frac{1}{2}det[g(v_i,v_j)].
$$
The Levi-Civita connection of $(M,g)$ determines a connection on the
bundle $\Lambda^2TM$, both denoted by $\nabla$, and the
corresponding curvatures are related by
$$
R(X\wedge Y)(Z\wedge T)=R(X,Y)Z\wedge T+Z\wedge R(X,Y)T
$$
for $X,Y,Z,T\in TM$. Let us note that we adopt the following definition for the curvature tensor
$R$ : $R(X,Y)=\nabla_{[X,Y]}-[\nabla_{X},\nabla_{Y}]$.

 The Hodge star operator defines an endomorphism $\ast$ of
$\Lambda^2TM$ with $\ast^2=Id$. Hence we have the decomposition
$$
\Lambda^2TM=\Lambda^2_{-}TM\oplus\Lambda^2_{+}TM
$$
where $\Lambda^2_{\pm}TM$ are the subbundles of $\Lambda^2TM$
corresponding to the $(\pm 1)$-eigenvalues of the operator $\ast$.

\smallskip

Let $(E_1,E_2,E_3,E_4)$ be a local oriented orthonormal frame of
$TM$. Set
\begin{equation}\label{s-basis}
s_1^{\pm}=E_1\wedge E_2\pm E_3\wedge E_4, \quad s_2^{\pm}=E_1\wedge E_3\pm E_4\wedge E_2, \quad
s_3^{\pm}=E_1\wedge E_4\pm E_2\wedge E_3.
\end{equation}
Then $(s_1^{\pm},s_2^{\pm},s_3^{\pm})$ is a local orthonormal frame of $\Lambda^2_{\pm}TM$ defining
an orientation on $\Lambda^2_{\pm}TM$,  which does not depend on the choice of the frame
$(E_1,E_2,E_3,E_4)$.

\smallskip

For every $a\in\Lambda ^2TM$, define a skew-symmetric endomorphism
of $T_{\pi(a)}M$ by

\begin{equation}\label{cs}
g(K_{a}X,Y)=2g(a, X\wedge Y), \quad X,Y\in T_{\pi(a)}M.
\end{equation}
Note that, denoting by $G$ the standard metric $-\frac{1}{2}Trace\,PQ$ on the space of
skew-symmetric endomorphisms, we have $G(K_a,K_b)=g(a,b)$ for $a,b\in \Lambda ^2TM$. If
$\sigma\in\Lambda^2_{+}TM$ is a unit vector, then $K_{\sigma}$ is a complex structure on the vector
space $T_{\pi(\sigma)}M$ compatible with the metric and the orientation of $M$. Conversely, the
$2$-vector $\sigma$ dual to one half of the K\"ahler $2$-form of such a complex structure is a unit
vector in $\Lambda^2_{+}TM$. Thus the unit sphere subbunlde ${\cal Z}_+={\cal Z}_+(M)$ of
$\Lambda^2_{+}TM$ parametrizes the complex structures on the tangent spaces of $M$ compatible with
its metric and orientation. This subbundle is called the twistor space of $M$.

\smallskip

The Levi-Civita connection $\nabla$ of $M$ preserves the bundles $\Lambda^2_{\pm}TM$, so it induces
a metric connection on these bundles denoted again by $\nabla$. The  horizontal distribution of
$\Lambda^2_{+}TM$ with respect to $\nabla$ is tangent to the twistor space ${\cal Z}_+$. Thus we
have the decomposition $T{\cal Z}_+={\cal H}\oplus {\cal V}$ of the tangent bundle of ${\cal Z}_+$
into horizontal and vertical components. The vertical space ${\cal V}_{\tau}=\{V\in T_{\tau}{\cal
Z}_+:~ \pi_{\ast}V=0\}$ at a point $\tau\in{\cal Z}_{+}$ is the tangent space to the fibre of ${\cal
Z}_+$ through $\tau$. Thus, considering $T_{\tau}{\cal Z}_+$ as a subspace of
$T_{\tau}(\Lambda^2_{+}TM)$ (as we shall always do), ${\cal V}_{\tau}$ is the orthogonal complement
of ${\Bbb R}\tau$ in $\Lambda^2_{+}T_{\pi(\tau)}M$. The map $V\ni{\cal V}_{\tau}\to K_{V}$ gives an
identification of the vertical space with the space of skew-symmetric endomorphisms of
$T_{\pi(\tau)}M$ that anti-commute with $K_{\tau}$. Let $s$ be a local section of ${\cal Z}_+$ such
that $s(p)=\tau$ where $p=\pi(\tau)$. Considering $s$ as a section of $\Lambda^2_{+}TM$, we have
$\nabla_{X}s\in{\cal V}_{\tau}$ for every $X\in T_pM$ since $s$ has a constant length. Moreover,
$X^h_{\tau}=s_{\ast}X-\nabla_{X}s$ is the horizontal lift of $X$ at ${\tau}$.

\smallskip

Denote by $\times$ the usual vector cross product on the oriented $3$-dimensional vector space
$\Lambda^2_{+}T_pM$, $p\in M$, endowed with the metric $g$. Then it is easy to check that
\begin{equation}\label{r-r}
g(R(a)b,c)=g({\cal R}(b\times c),a)
\end{equation}
for $a\in\Lambda^2T_pM$, $b,c\in\Lambda^2_{+}T_pM$. It is also easy to show that for every
$a,b\in\Lambda^2_{+}T_pM$
\begin{equation}\label{com}
K_a\circ K_b=-g(a,b)Id+ K_{a\times b}.
\end{equation}

For every $t>0$, define a Riemannian metric $h_t$ by
$$
h_t(X^h_{\sigma}+V,Y^h_{\sigma}+W)=g(X,Y)+tg(V,W)
$$
for $\sigma\in{\cal Z}_{+}$, $X,Y\in T_{\pi(\sigma)}M$, $V,W\in{\cal V}_{\sigma}$.
\smallskip

The twistor space ${\cal Z}_+$ admits two natural almost complex structures that are compatible
with the metrics $h_t$. One of them has been introduced by Atiyah, Hitchin and Singer who have
proved that it is integrable if and only if the base manifold is anti-self-dual \cite{AHS}. The
other one, introduced by Eells and Salamon, although never integrable, plays an important role in
harmonic maps theory \cite{ES}.

\smallskip

 The action of $SO(4)$ on $\Lambda^2{\Bbb R}^4$ preserves the
decomposition $\Lambda^2{\Bbb R}^4=\Lambda^2_{+}{\Bbb R}^4\oplus \Lambda^2_{-}{\Bbb R}^4$. Thus,
considering $S^2$ as the unit sphere in $\Lambda^2_{+}{\Bbb R}^4$, we have an action of the group
$SO(4)$ on $S^2$. Then, if $SO(M)$ denotes  the principal bundle of the oriented orthonormal frames
on $M$, the twistor space ${\cal Z}_+={\cal Z}_+(M)$ is the associated bundle $SO(M)\times_{SO(4)}
S^2$. It follows from the Vilms theorem (see, for example, \cite[Theorem 9.59]{Besse}) that the
projection map $\pi:({\cal Z}_+,h_t)\to (M,g)$ is a Riemannian submersion with totally geodesic
fibres (this can also be proved by a direct computation).

\smallskip

Denote by $D$ the Levi-Chivita connection of $({\cal Z}_+,h_t)$.

\smallskip

Let $(N,x_1,...,x_4)$ be a local coordinate system of $M$ and let $(E_1,...,E_4)$ be an oriented
orthonormal frame of $TM$ on $N$. If $(s_1^+,s_2^+,s_3^+)$ is the local frame of $\Lambda^2_{+}TM$
define by (\ref{s-basis}), then $\widetilde x_{a}=x_{a}\circ\pi$, $y_j(\tau)=g(\tau,
(s_j^+\circ\pi)(\tau))$, $1\leq a \leq 4$, $1\leq j\leq 3$, are local coordinates of
$\Lambda^2_{+}TM$ on $\pi^{-1}(N)$.

   The horizontal lift $X^h$ on $\pi^{-1}(N)$ of a vector field
$$
X=\sum_{a=1}^4 X^{a}\frac{\partial}{\partial x_{a}}
$$
is given by
\begin{equation}\label{hl}
X^h=\sum_{a=1}^4 (X^{a}\circ\pi)\frac{\partial}{\partial \widetilde{x}_{a}}
-\sum_{j,k=1}^3y_j(g(\nabla_{X}s_j,s_k)\circ\pi)\frac{\partial}{\partial y_k}.
\end{equation}
Hence
\begin{equation}\label{Lie-1}
[X^h,Y^h]=[X,Y]^h+\sum_{j,k=1}^3y_j(g(R(X\wedge Y)s_j,s_k)\circ\pi)\frac{\partial}{\partial y_k}
\end{equation}
for every vector fields $X,Y$ on $N$. Let $\tau\in{\cal Z}_+$. Using the standard identification
$T_{\tau}(\Lambda^2_{+}T_pM)\cong \Lambda^2_{+}T_{\pi(\tau)}M$  we obtain from (\ref{Lie-1}) the
well-known formula
\begin{equation}\label{Lie-2}
[X^h,Y^h]_{\tau}=[X,Y]^h_{\tau}+R_{p}(X\wedge Y)\tau, \quad p=\pi(\tau).
\end{equation}

Then we have the following
\begin{lemma}\label{LC} {\rm (\cite{DM})}
If $X,Y$ are (local) vector fields on $M$ and $V$ is a vertical vector field on ${\cal Z}_+$, then
\begin{equation}\label{D-hh}
(D_{X^h}Y^h)_{\tau}=(\nabla_{X}Y)^h_{\tau}+\frac{1}{2}R_{p}(X\wedge Y)\tau,
\end{equation}
\begin{equation}\label{D-vh}
(D_{V}X^h)_{\tau}={\cal H}(D_{X^h}V)_{\tau}=-\frac{t}{2}(R_{p}(\tau\times V)X)^h_{\tau}
\end{equation}
where $\tau\in{\cal Z}_+$, $p=\pi(\tau)$, and ${\cal H}$ means "the horizontal component".
\end{lemma}
{\bf Proof}. Identity (\ref{D-hh}) follows from the Koszul formula for the Levi-Chivita connection
and (\ref{Lie-2}).

Let $W$ be a vertical vector field on ${\cal Z}_+$. Then
$$
h_t(D_{V}X^h,W)=-h_t(X^h,D_{V}W)=0
$$
since the fibres are totally geodesic submanifolds, so $D_{V}W$ is a vertical vector field.
Therefore $D_{V}X^h$ is a horizontal vector field. Moreover, $[V,X^h]$ is a vertical vector field,
hence $D_{V}X^h={\cal H}D_{X^h}V$. Thus
$$
h_t(D_{V}X^h,Y^h)=h_t(D_{X^h}V,Y^h)=-h_t(V,D_{X^h}Y^h).
$$
Now (\ref{D-vh}) follows from (\ref{D-hh}) and (\ref{r-r}).

\section{A hypersurface in ${\cal Z}_+$ determined by an almost Hermitian structure on $M$}

Let $(g,J)$ be an almost Hermitian structure on a four-manifold $M$. Define a section $\alpha$ of
$\Lambda^2TM$ by
$$
g(\alpha,X\wedge Y)=\frac{1}{2}g(JX,Y),~X,Y\in TM.
$$
Thus, at any point of $M$, $\alpha$ is the dual $2$-vector of one half of the K\"ahler $2$-form of
the almost Hermitian manifold $(M,g,J)$. Note also that $K_{\alpha_p}=J_p$ for every $p\in M$.

Consider $M$ with the orientation yielded by the almost complex structure $J$. Then $\alpha$ is a
section of the twistor bundle ${\cal Z}_+$. This section determines a hypersurface of the twistor
space defined by
$$
\Sigma_J=\{\sigma\in{\cal Z}_+:~g(\sigma,\alpha_{\pi(\sigma)})=0\}.
$$
By (\ref{com}), the points of $\Sigma_J$ are complex structures on the tangent spaces of $M$ that
are compatible with the metric and the orientation, and anti-commute with $J$.

 Clearly, $\Sigma_J$ is the circle bundle of the rank 2 vector bundle
$$\Lambda^2_{0}=\{\sigma\in\Lambda^2_{+}TM:~g(\sigma,\alpha_{\pi(\sigma)})=0\}.$$
As is well-known (and easy to see),  the complexification of this bundle is the bundle
$\Lambda^{2,0}\oplus\Lambda^{0,2}$ where $\Lambda^{r,s}$ stands for the bundle of $(r+s)$-vectors
of type $(r,s)$ with respect to $J$.

We shall compute the second fundamental form $\Pi$ of the hypersurface $\Sigma_J$ in $({\cal
Z}_+,h_t)$.

Note that for $\sigma\in\Sigma_J$
$$
T_{\sigma}\Sigma_J=\{E\in T_{\sigma}{\cal Z}_+:~g({\cal
V}E,\alpha_{\pi(\sigma)})=-g(\sigma,\nabla_{\pi_{\ast E}}\alpha)\}
$$
where ${\cal V}E$ means "the vertical component of $E$". Therefore
$$
T_{\sigma}\Sigma_J=\{X^h_{\sigma}-g(\sigma,\nabla_{X}\alpha)\alpha_{\pi(\sigma)}:~
X\in T_{\pi(\sigma)}M\}\oplus{\Bbb
R}(\alpha_{\pi(\sigma)}\times\sigma).
$$

 Given $\tau\in{\cal Z}$ and $X\in T_{\pi(\tau)}M$, define a vertical vector of ${\cal Z}_+$ at $\tau$
by
$$
X^v_{\tau}=-g(\tau,\nabla_{X}\alpha)\alpha_{\pi(\tau)}+g(\tau,\alpha_{\pi(\tau)})\nabla_{X}\alpha.
$$

Set
$$
\widehat X_{\tau}=X^h_{\tau}+X^v_{\tau}.
$$
Thus every (local) vector field $X$ on $M$, gives rise to a vector field $\widehat X$ on ${\cal
Z}_+$ tangent to $\Sigma_J$.

\smallskip

Let $\rho(\tau)=g(\tau,\alpha_{\pi(\tau)})$, $\tau\in{\cal Z}_+$,  be the defining function of
$\Sigma_J$ and let $grad\,\rho$ be the gradient vector field of the function $\rho$ with respect to
the metric $h_t$. Fix a point $\tau\in{\cal Z}_+$ and take a section $s$ of ${\cal Z}_+$ such that
$s_{\pi(\tau)}=\tau$,  $\nabla s|_{\pi(\tau)}=0$. Then, for $X\in T_{\pi(\tau)}M$,
\begin{equation}\label{gr-h}
h_t(X^h_{\tau},grad\,\rho)=s_{\ast}(X)(\rho)=X(g(s,\alpha))=g(\tau,\nabla_X\alpha).
\end{equation}
Moreover, if $V\in{\cal V}_{\tau}$,
\begin{equation}\label{gr-v}
h_t(V,grad\,\rho)=V(\sum_{k=1}^3y_k(g(s_k,\alpha)\circ\pi))=\sum_{k=1}^3V(y_k)g(s_k,\alpha)_{\pi(\tau)}
=g(V,\alpha_{\pi(\tau)}).
\end{equation}

\begin{lemma}\label{hat-hat}
If $\sigma\in\Sigma_J$ and $X,Y\in T_{\pi(\sigma)}M$, then
$$
\begin{array}{c}
h_t(\Pi(\widehat X,\widehat Y),grad\,\rho)_{\sigma}=
\displaystyle{\frac{t}{2}}[g(\sigma,\nabla_{X}\alpha)g(\sigma,\nabla_{R(\sigma\times\alpha_{\pi(\sigma)})Y}\alpha)\\[10pt]
\hspace{4cm}+g(\sigma,\nabla_{Y}\alpha)g(\sigma,\nabla_{R(\sigma\times\alpha_{\pi(\sigma)})X}\alpha)]\\[10pt]
\hspace{4cm}\displaystyle{-\frac{1}{2}g(\sigma,\nabla^2_{XY}\alpha)-\frac{1}{2}g(\sigma,\nabla^2_{YX}\alpha)}
\end{array}
$$
where $\nabla^{2}_{XY}\alpha=\nabla_X\nabla_Y\alpha-\nabla_{\nabla_XY}\alpha$ is the second
covariant derivative of $\alpha$.
\end{lemma}
{\bf Proof}. Extend $X$ and $ Y$ to vector fields in a neighbourhood of the point $p=\pi(\sigma)$.
It follows from (\ref{D-hh}), (\ref{gr-h}) and (\ref{gr-v}) that
\begin{equation}\label{hhgr}
h_t(D_{X^h}Y^h,grad\,\rho)_{\sigma}=g(\nabla_{\nabla_XY}\alpha,\sigma)+\frac{1}{2}g(R(X\wedge
Y)\sigma,\alpha_{p}).
\end{equation}
Identities (\ref{D-vh}) and (\ref{gr-h}) imply
\begin{equation}\label{vhgr}
\begin{array}{c}
h_t(D_{X^v}Y^h,grad\,\rho)_{\sigma}=
\displaystyle{\frac{t}{2}}g(\sigma,\nabla_X\alpha)g(\sigma,\nabla_{R(\sigma\times\alpha_p)Y}\alpha).
\end{array}
\end{equation}

Next, note that
$$
h_t(D_{X^h}Y^v,grad\,\rho)=h_t([X^h,Y^v],grad\,\rho)+h_t(D_{Y^v}X^h,grad\,\rho).
$$
Take an oriented orthonormal frame $(E_1,...,E_4)$ of $M$ near $p$ such that $\nabla E_{a}|_{p}=0$,
$a=1,...,4$. Then $\nabla s_i^+|_{p}=0$, $i=1,2,3$, which implies
\begin{equation}\label{L}
X^h_{\sigma}=\sum_{a=1}^4X^{a}(p)\frac{\partial}{\partial \widetilde{x}_{a}}(\sigma),\quad
[X^h,\frac{\partial}{\partial y_i}]_{\sigma}=0,\>i=1,2,3.
\end{equation}
We have
\begin{equation}\label{vl}
Y^v=\sum_{j,k=1}^3
y_k(g(s_k^+,\alpha)g(\nabla_{Y}\alpha,s_j^+)-g(s_j^+,\alpha)g(\nabla_{Y}\alpha,s_k^+))\circ\pi
\frac{\partial }{\partial y_j}.
\end{equation}
It follows from (\ref{L}) and (\ref{vl}) that
$$
[X^h,Y^v]_{\sigma}=g(\sigma,\nabla_X\alpha)\nabla_Y\alpha
-g(\sigma,\nabla_Y\alpha)\nabla_X\alpha-g(\sigma,\nabla_X\nabla_Y\alpha)\alpha_p.
$$
Hence, by (\ref{gr-v}),
$$
h_t([X^h,Y^v],grad\,\rho)_{\sigma}=-g(\sigma,\nabla_X\nabla_Y\alpha).
$$
Thus we have
\begin{equation}\label{hvgr}
h_t(D_{X^h}Y^v,grad\,\rho)_{\sigma}=-g(\sigma,\nabla_X\nabla_Y\alpha)+
\displaystyle{\frac{t}{2}}g(\sigma,\nabla_Y\alpha)g(\sigma,\nabla_{R(\sigma\times\alpha_p)X}\alpha).
\end{equation}

 The fibres of of ${\cal Z}_+$ are totally geodesic
submanifolds, hence $(D_{X^v}Y^v)_{\sigma}$ is the standard covariant derivative on the unit sphere
in the vector space $\Lambda^2_{+}T_{p}M$. It follows from (\ref{vl}) that
$$
\begin{array}{c}
(D_{X^v}Y^v)_{\sigma}=g(X^v_{\sigma},\alpha_{p})[\nabla_Y\alpha-g(\nabla_Y\alpha,\sigma)\sigma]
-g(X^v_{\sigma},\nabla_Y\alpha)[\alpha_{p}-g(\alpha_{p},\sigma)\sigma]=\\[10pt]
-g(\sigma,\nabla_X\alpha)[\nabla_Y\alpha-g(\nabla_Y\alpha,\sigma)\sigma].
\end{array}
$$
Hence
\begin{equation}\label{vvgr}
h_t(D_{X^v}Y^v,grad\,\rho)_{\sigma}=0.
\end{equation}
Now the lemma follows from identities (\ref{hhgr}),  (\ref{vhgr}), (\ref{hvgr}), and (\ref{vvgr}).

\smallskip

If $\sigma\in\Sigma_J$, the vertical part of $T_{\sigma}\Sigma_J$ is ${\Bbb
R}(\alpha_{\pi(\sigma)}\times\sigma)$. Define a vertical vector field $\xi$ on ${\cal Z}_+$ tangent
to $\Sigma_J$ setting
$$
\xi_{\tau}=\alpha_{\pi(\tau)}\times\tau,\quad \tau\in{\cal Z}_+.
$$
\begin{lemma}\label{hat-x}
If $\sigma\in\Sigma_J$ and $X\in T_{\pi(\sigma)}M$, then
$$
h_t(\Pi(\xi,\widehat
X),grad\,\rho)_{\sigma}=-g(\xi_{\sigma},\nabla_X\alpha)-\frac{t}{2}g(\sigma,\nabla_{R(\alpha_{\pi(\sigma)})X}\alpha).
$$
$$
h_t(\Pi(\xi,\xi),grad\,\rho)_{\sigma}=0.
$$
\end{lemma}
{\bf Proof}. Identity (\ref{D-vh}) implies
\begin{equation}\label{xhgr}
h_t(D_{\xi}X^h,grad\,\rho)_{\sigma}=-\displaystyle{\frac{t}{2}}g(\sigma,\nabla_{R(\sigma\times\xi_{\sigma})X}\alpha).
\end{equation}
A simple computation gives
$$
(D_{\xi}X^v)_{\sigma}=-g(\xi_{\sigma},\nabla_X\alpha)\alpha_{\pi(\sigma)},\quad
(D_{\xi}\xi)_{\sigma}=0.
$$
Hence
\begin{equation}\label{xi-gr}
h_t(D_{\xi}X^v,grad\,\rho)_{\sigma}=-g(\xi_{\sigma},\nabla_X\alpha),\quad
h_t(D_{\xi}\xi,grad\,\rho)_{\sigma}=0.
\end{equation}
Thus the result follows from (\ref{xhgr}) and (\ref{xi-gr}).

\begin{prop}\label{Pi}
Let $\sigma\in\Sigma_J$ and $E,F\in T_{\sigma}\Sigma_J$. Set $X=\pi_{\ast}E$, $Y=\pi_{\ast}F$,
$V={\cal V}E$, $W={\cal V}F$. Then
$$
\begin{array}{c}
h_t(\Pi(E,F),grad\,\rho)_{\sigma}=\\[6pt]
\displaystyle{\frac{t}{2}g(\sigma,\nabla_{X}\alpha)g(\sigma,\nabla_{R(\sigma\times\alpha_{\pi(\sigma)})Y}\alpha)
+\frac{t}{2}g(\sigma,\nabla_{Y}\alpha)g(\sigma,\nabla_{R(\sigma\times\alpha_{\pi(\sigma)})X}\alpha)}\\[8pt]
\displaystyle{-\frac{1}{2}g(\sigma,\nabla^2_{XY}\alpha)-\frac{1}{2}g(\sigma,\nabla^2_{YX}\alpha)}\\[8pt]
+\displaystyle{\frac{t}{2}g(\alpha_{\pi(\sigma)}\times
V,\nabla_{R(\alpha_{\pi(\sigma)})Y}\alpha)+\frac{t}{2}g(\alpha_{\pi(\sigma)}\times W,\nabla_{R(\alpha_{\pi(\sigma)})X}\alpha)}\\[8pt]
-g(V,\nabla_Y\alpha)-g(W,\nabla_X\alpha).
\end{array}
$$
\end{prop}
{\bf Proof}. This follows from Lemmas~\ref{hat-hat} and \ref{hat-x} taking into account that
$E=\widehat X_{\sigma}+g(V,\xi_{\sigma})\xi_{\sigma}$, $F=\widehat
Y_{\sigma}+g(W,\xi_{\sigma})\xi_{\sigma}$.

\begin{cor}
If $(M,g,J)$ is K\"ahler, $\Sigma_J$ is a totally geodesic submanifold of ${\cal Z}_{+}$.
\end{cor}

\section{Minimality of the hypersurface $\Sigma_J$}

Let $\Omega(X,Y)=g(JX,Y)$ be the fundamental $2$-form of the almost Hermitian manifold $(M,g,J)$.
Denote by $N$ the Nijenhuis tensor of $J$, $N(Y,Z)=-[Y,Z]+[JY,JZ]-J[Y,JZ]-J[JY,Z]$. It is
well-known (and easy to check) that
\begin{equation}\label{nJ}
2g((\nabla_XJ)(Y),Z)=d\Omega(X,Y,Z)-d\Omega(X,JY,JZ)+g(N(Y,Z),JX).
\end{equation}

\subsection{The case of integrable $J$}

Suppose that the almost complex structure $J$ is integrable. Note that the integrability condition
for $J$ is equivalent to $(\nabla_{X}J)(Y)=(\nabla_{JX}J)(JY)$, $X,Y\in TM$ \cite[Corollary
4.2]{G}. Let $B$ be the vector field on $M$ dual to the Lee form $\theta=-\delta\Omega\circ J$ with
respect to the metric $g$. Then (\ref{nJ}) and the identity $d\Omega=\Omega\wedge\theta$ imply the
following well-known formula
\begin{equation}\label{nablaJ}
2(\nabla_XJ)(Y)=g(JX,Y)B-g(B,Y)JX+g(X,Y)JB-g(JB,Y)X.
\end{equation}
We have $g(\nabla_X\alpha,Y\wedge Z)=\displaystyle{\frac{1}{2}}g((\nabla_XJ)(Y),Z)$ and it follows
that
\begin{equation}\label{nalpha}
\nabla_X\alpha=\displaystyle{\frac{1}{2}}(JX\wedge B+X\wedge JB).
\end{equation}
The latter identity implies
\begin{equation}\label{n2alpha}
\nabla^2_{XY}\alpha=\displaystyle{\frac{1}{2}}[(\nabla_XJ)(Y)\wedge
B+Y\wedge(\nabla_XJ)(B)+JY\wedge\nabla_XB+Y\wedge J\nabla_XB].
\end{equation}
Let $\sigma\in\Sigma_J$ and $X,Y\in T_{\pi(\sigma})M$. Then a simple computation using identities
(\ref{cs}), (\ref{com}), (\ref{nablaJ}) - (\ref{n2alpha}) gives
$$
\begin{array}{c}
g(\sigma,\nabla_X\alpha)=\displaystyle{\frac{1}{2}}g(X,K_{\xi_{\sigma}}B),\\[6pt]
4g(\sigma,\nabla^2_{XY}\alpha)=-g(K_{\xi_{\sigma}}B\wedge B,X\wedge Y+JX\wedge
JY)-g(JX,B)g(JY,K_{\xi_{\sigma}}B)\\[6pt]
+\displaystyle{\frac{1}{2}}||B||^2g(X,K_{\xi_{\sigma}}Y)-2g(\nabla_XB,K_{\xi_{\sigma}}Y)
\end{array}
$$
where, as above, $\xi_{\sigma}=\alpha_{\pi(\sigma)}\times\sigma$. Moreover, if $V\in
T_{\sigma}\Sigma_J$ is a vertical vector,
$$
g(\alpha_{\pi(\sigma)}\times  V,\nabla_X\alpha)=-\frac{1}{2}g(V,\xi_{\sigma})g(\sigma,JX\wedge
B+X\wedge JB)=-\frac{1}{2}g(V,\xi_{\sigma})g(X,K_{\xi_{\sigma}}B),
$$
$$
g(V,\nabla_X\alpha)=-\frac{1}{2}g(V,\xi_{\sigma})g(X,K_{\sigma}B).
$$
Now Proposition~\ref{Pi} can be rewritten as
\begin{prop}\label{Pi-new}
Suppose that the almost complex structure $J$ is integrable. Let $\sigma\in\Sigma_J$ and $E,F\in
T_{\sigma}\Sigma_J$. Set $X=\pi_{\ast}E$, $Y=\pi_{\ast}F$, $V={\cal V}E$, $W={\cal V}F$ and
$\xi_{\sigma}=\alpha_{\pi(\sigma)}\times\sigma$. Then
$$
\begin{array}{c}
h_t(\Pi(E,F),grad\,\rho)_{\sigma}=\\[6pt]
-\displaystyle{\frac{t}{8}g(X,K_{\xi_{\sigma}}B)g(Y,R(\xi_{\sigma})K_{\xi_{\sigma}}B)-\frac{t}{8}g(Y,K_{\xi_{\sigma}}B)g(X,R(\xi_{\sigma})K_{\xi_{\sigma}}B)}\\[8pt]
+\displaystyle{\frac{1}{8}g(JX,B)g(JY,K_{\xi_{\sigma}}B)+\frac{1}{8}g(JY,B)g(JX,K_{\xi_{\sigma}}B)}\\[8pt]
+\displaystyle{\frac{1}{4}g(\nabla_XB,K_{\xi_{\sigma}}Y) +\frac{1}{4}g(\nabla_YB,K_{\xi_{\sigma}}X)}\\[8pt]
\displaystyle{-\frac{1}{4}g(V,\xi_{\sigma})g(R(\alpha_{\pi(\sigma)})Y,K_{\xi_{\sigma}}B)-\frac{1}{4}g(W,\xi_{\sigma})g(R(\alpha_{\pi(\sigma)})X,K_{\xi_{\sigma}}B)}\\[8pt]
+\displaystyle{\frac{1}{4}g(V,\xi_{\sigma})g(X,K_{\sigma}B)+\frac{1}{4}g(W,\xi_{\sigma})g(Y,K_{\sigma}B)}
\end{array}
$$
\end{prop}
\begin{cor}\label{trace}
Let $\sigma\in\Sigma_J$. Then
$$
h_t(Trace\,\Pi,grad\,\rho)_\sigma=\frac{1}{2}(d\theta+\theta\wedge
d\ln\sqrt{8+2t||\theta_{\pi(\sigma)}||^2})(\alpha_{\pi(\sigma)}\times\sigma).
$$
\end{cor}
{\bf Proof}. Set $p=\pi(\sigma)$. Suppose first that $B_p\neq 0$. Then $E_1=||B_p||^{-1}B_p$,
$E_2=K_{\alpha_{p}}E_1$, $E_3=K_{\sigma}E_1$, $E_4=K_{\xi_{\sigma}}E_1$ form an oriented
orthonormal basis of $T_{p}M$ such that $\alpha_{p}=(s_1)^{+}_{p}$, $\sigma=(s_2)^{+}_{p}$,
$\xi_{\sigma}=(s_3)^{+}_{p}$ where $s_1^+,s_2^+,s_3^+$ is the basis of $\Lambda^2_{+}T_{p}M$
defined by means of $E_1,...,E_4$ via (\ref{s-basis}). We have
$$
g(\sigma,\nabla_X\alpha)=\frac{1}{2}g(X,K_{\xi_{\sigma}}B)=\frac{1}{2}||B_p||g(X,E_4),\quad X\in
T_pM.
$$
Hence $\widehat{(E_i)}_{\sigma}=(E_i)^h_{\sigma}$ for $i=1,2,3$ and
$\widehat{(E_4)}_{\sigma}=(E_4)^h_{\sigma}-\displaystyle{\frac{1}{2}}||B_p||\alpha_{p}$. Thus
$\widehat{(E_i)}_{\sigma}$, $i=1,2,3$,
$(1+\displaystyle{\frac{t}{4}}||B_p||^2)^{-\frac{1}{2}}\widehat{(E_4)}_{\sigma}$,
$\displaystyle{\frac{1}{\sqrt t}}\xi_{\sigma}$ constitute an orthonormal basis of
$T_{\sigma}\Sigma_J$. Note that
$$
g(E_4,R(\xi_{\sigma})K_{\xi_{\sigma}}B)=||B_p||^{-1}g(K_{\xi_{\sigma}}B,R(\xi_{\sigma})K_{\xi_{\sigma}}B)=0
$$
and
$$
g(JE_4,B)=||B_p||^{-1}g(K_{\alpha_p}\circ K_{\xi_{\sigma}}B,B)=-||B_p||^{-1}g(K_{\sigma}B,B)=0.
$$
Then, by Proposition~\ref{Pi-new},
$$
\begin{array}{c}
h_t(Trace\,\Pi,grad\,\rho)_\sigma=
-\displaystyle{\frac{1}{4}g(K_{\xi_{\sigma}}B,R(\xi_{\sigma})K_{\xi_{\sigma}}B)+\frac{1}{4}g(B,K_{\xi_{\sigma}}B)}\\[6pt]
+\displaystyle{\frac{1}{2}\sum_{i=1}^3g(\nabla_{E_i}B,K_{\xi_{\sigma}}E_i)+\frac{1}{2}(1+\frac{t}{4}||B_p||^2)^{-1}g(\nabla_{E_4}B,K_{\xi_{\sigma}}E_4)}\\[6pt]
=\displaystyle{\frac{1}{2}\sum_{j=1}^4g(\nabla_{E_j}B,K_{\xi_{\sigma}}E_j)+\frac{t||B_p||}{8+2t||B_p||^2}g(\nabla_{E_4}B,B)}\\[6pt]
=\displaystyle{\frac{1}{2}d\theta(E_1\wedge E_4+E_2\wedge
E_3)+\frac{t}{16+4t||B_p||^2}(\theta\wedge d||\theta||^2)(E_1\wedge E_4+E_2\wedge E_3)}
\end{array}
$$
If $B_p=0$, then $\nabla_{X}\alpha=0$ for every $X\in T_pM$ by (\ref{nalpha}). Taking a unit vector
$E_1\in T_pM$ we set $E_2=K_{\alpha_{p}}E_1$, $E_3=K_{\sigma}E_1$, $E_4=K_{\xi_{\sigma}}E_1$. Then
$\widehat{(E_j)}_{\sigma}=(E_j)^h_{\sigma}$, $j=1,...,4$, $\displaystyle{\frac{1}{\sqrt
t}}\xi_{\sigma}$ constitute an orthonormal basis of $T_{\sigma}\Sigma_J$. It follows from
Proposition~\ref{Pi-new} that
$$
h_t(Trace\,\Pi,grad\,\rho)_\sigma=
\displaystyle{\frac{1}{2}\sum_{j=1}^4g(\nabla_{E_j}B,K_{\xi_{\sigma}}E_j)=\frac{1}{2}d\theta(E_1\wedge
E_4+E_2\wedge E_3)}.
$$
\begin{prop}\label{min}
If $J$ is integrable, the hypersurface $\Sigma_J$ is a minimal submanifold of $({\cal Z}_+,h_t)$ if
and only if the $2$-form
$$
d\frac{\theta}{\sqrt{8+2t||\theta||^2}}
$$
is of type $(1,1)$ with respect to $J$.
\end{prop}
{\bf Proof}. The condition that $\Sigma_J$ is a minimal submanifold means that\\
$h_t(Trace\,\Pi,grad\,\rho)=0$ on $\Sigma_J$.

 Let $p\in M$ and take an orthonormal basis of $T_pM$ of
the form $E_1, E_2=JE_1, E_3, E_4=JE_3$. Then $\alpha_p=(s_1)^{+}_p$, so $(s_2)^{+}_p$,
$(s_3)^{+}_p\in\Sigma_J$.

It is easy to check that, for every $a\in\Lambda^2_{+}T_pM$ and $b\in\Lambda^2_{-}T_pM$, the
endomorphisms $K_a$ and $K_b$ of $T_pM$ commute. It follows that, for every $X,Y\in T_pM$, the
$2$-vector $X\wedge Y-JX\wedge JY$ is orthogonal to $\Lambda^2_{-}T_pM$, therefore it lies in
$\Lambda^2_{+}T_pM$. Moreover, $X\wedge Y-JX\wedge JY$ is orthogonal to $\alpha_p$, hence is a
linear combination of $(s_2)^{+}_p=-\alpha_p\times (s_3)^{+}_p$ and $(s_3)^{+}_p=\alpha_p\times
(s_2)^{+}_p$. Thus if $h_t(Trace\,\Pi,grad\,\rho)=0$ on $\Sigma_J$,
$$
(d\theta+\theta\wedge d\ln\sqrt{8+2t||\theta||^2})(X\wedge Y-JX\wedge JY)=0.
$$
Conversely, if this identity holds, then $h_t(Trace\,\Pi,grad\,\rho)=0$ at the points $(s_2)^{+}_p$
and $(s_3)^{+}_p$ of $\Sigma_J$. For every $\sigma\in\Sigma_J$ with $\pi(\sigma)=p$, the $2$-vector
$\alpha_p\times\sigma$ is a linear combination of $(s_2)^{+}_p$ and $(s_3)^{+}_p$ , hence
$h_t(Trace\,\Pi,grad\,\rho)=0$ on $\Sigma_J$.

Thus, $\Sigma_J$ is minimal if and only if the form $d\theta+\theta\wedge
d\ln\sqrt{8+2t||\theta||^2}$ is of type $(1,1)$. But the this condition is equivalent to the
condition that the form $d\displaystyle{\frac{\theta}{\sqrt{8+2t||\theta||^2}}}$ is of type
$(1,1)$.

\medskip

\subsection{The case of symplectic $J$} Now suppose that $d\Omega=0$. Then, by (\ref{nJ}),
$$
g((\nabla_XJ)(Y),Z)=\frac{1}{2}g(N(Y,Z),JX).
$$
The Nijenhuis tensor $N(Y,Z)$ is skew-symmetric, so it induces a linear map $\Lambda^2TM\to TM$
which we denote again by $N$. The identity
\begin{equation}\label{nalpha-sympl}
g(\nabla_X\alpha,Y\wedge Z)=\frac{1}{4}g(N(Y,Z),JX)
\end{equation}
implies that, for every $a\in\Lambda^2TM$ and $X\in T_{\pi(a)}M$,
$$
g(\nabla_X\alpha, a)=\frac{1}{4}g(N(a),JX).
$$
In particular, if $V\in{\cal V}_{\sigma}$, then
$$
g(\nabla_X\alpha, V)=\frac{1}{4}g(V,\xi_{\sigma})g(N(\xi_{\sigma}),JX).
$$

 Proposition~\ref{Pi} implies the following.
\begin{prop}\label{Pi-sympl}
Suppose that $d\Omega=0$. Let $\sigma\in\Sigma_J$ and $E,F\in T_{\sigma}\Sigma_J$. Set
$X=\pi_{\ast}E$, $Y=\pi_{\ast}F$, $V={\cal V}E$, $W={\cal V}F$ and
$\xi_{\sigma}=\alpha_{\pi(\sigma)}\times\sigma$. Then
$$
\begin{array}{c}
h_t(\Pi(E,F),grad\,\rho)_{\sigma}=\\[6pt]
\displaystyle{\frac{t}{32}}g(JN(\sigma),X)g(JN(\sigma),R(\xi_{\sigma})Y)
+\displaystyle{\frac{t}{32}}g(JN(\sigma),Y)g(JN(\sigma),R(\xi_{\sigma})X)\\[8pt]
\displaystyle{-\frac{1}{2}}g(\sigma,\nabla^2_{XY}\alpha)-\frac{1}{2}g(\sigma,\nabla^2_{YX}\alpha)\\[8pt]
+\displaystyle{\frac{t}{8}g(N(\alpha_{\pi(\sigma)}\times V),JR(\alpha_{\pi(\sigma)})Y)
+\frac{t}{8}g(N(\alpha_{\pi(\sigma)}\times W),JR(\alpha_{\pi(\sigma)})X)}\\[8pt]
-\displaystyle{\frac{1}{4}g(V,\xi_{\sigma})g(N(\xi_{\sigma}),JX)-\frac{1}{4}g(W,\xi_{\sigma})g(N(\xi_{\sigma}),JY)}.
\end{array}
$$
\end{prop}
\begin{cor}\label{trace-sympl}
Let $\sigma\in\Sigma_J$. Then
$$
h_t(Trace\,\Pi,grad\,\rho)_\sigma=-g(trace\nabla^2\alpha,\sigma).
$$
\end{cor}
{\bf Proof}.  Suppose first that $N(\sigma)\neq 0$. Take an orthonormal basis of $T_{\pi(\sigma)}M$
of the form $E_1$, $E_2=JE_1$, $E_3=||N(\sigma)||^{-1}N(\sigma)$, $E_4=JE_3$. Then, by
(\ref{nalpha-sympl}),
$$
\widehat{(E_k)}_{\sigma}=(E_k)^h_{\sigma}+\frac{1}{4}||N(\sigma)||g(E_4,E_k)\alpha_{\pi(\sigma)},
$$
$k=1,...,4$. Thus $\widehat{(E_i)}_{\sigma}$, $i=1,2,3$,
$(1+\displaystyle{\frac{t}{16}}||N(\sigma)||^2)^{-\frac{1}{2}}\widehat{(E_4)}_{\sigma}$,
$\displaystyle{\frac{1}{\sqrt t}}\xi_{\sigma}$ form an orthonormal basis of $T_{\sigma}\Sigma_J$.
Note also that
$$
g(JN(\sigma),R(\xi_{\sigma})E_4)=||N(\sigma)||g(E_4,R(\xi_{\sigma})E_4)=0.
$$
Then Proposition~\ref{Pi-sympl} implies
$$
\begin{array}{c}
h_t(Trace\,\Pi,grad\,\rho)_\sigma=-\displaystyle{\frac{t}{16}}g(JN(\sigma),R(\xi_{\sigma})JN(\sigma))-g(trace\nabla^2\alpha,\sigma)\\[8pt]
=-g(trace\nabla^2\alpha,\sigma).
\end{array}
$$
If $N(\sigma)=0$, then, in view of (\ref{nalpha-sympl}),  $g(\nabla_X\alpha, \sigma)=0$. Thus
$\widehat{(E_k)}_{\sigma}=(E_k)^h_{\sigma}$ for any orthonormal basis $E_k$ of $T_{\pi(\sigma)}M$,
$k=1,...,4$, and the result is a direct consequence of Proposition~\ref{Pi-sympl}.

Denote by $\rho^{\ast}$ the $\ast$-Ricci tensor of the almost Hermitian manifold $(M,g,J)$. Recall
that it is defined as $\rho^{\ast}(X,Y)=trace\{Z\to R(JZ,X)JY\}$. Note that
\begin{equation}\label{rho-star}
\rho^{\ast}(JX,JY)=\rho^{\ast}(Y,X),
\end{equation}
in particular $\rho^{\ast}(X,JX)=0$.
\begin{prop}\label{min-sympl}
If $d\Omega=0$, then the hypersurface $\Sigma_J$ is a minimal submanifold of $({\cal Z}_+,h_t)$ if
and only if the tensor  $\rho^{\ast}$ is symmetric.
\end{prop}
{\bf Proof}. The form $\Omega$ is harmonic since $d\Omega=0$ and $\ast\Omega=\Omega$. Then, by
Corollary~\ref{trace-sympl} and the Weitzenb\"ock formula, $\Sigma_J$ is  minimal if and only if,
for every $2$-form $\tau\in\Lambda^2_{+}T^{\ast}M$ orthogonal to $\Omega$, $g(S(\Omega),\tau)=0$
where
$$
S(\Omega)(X,Y)=Trace\{Z\to (R(Z,Y)\Omega)(Z,X)-(R(Z,X)\Omega)(Z,Y)\}
$$
(see, for example, \cite{EL}). We have
$$
\begin{array}{c}
(R(Z,Y)\Omega)(Z,X)=-\Omega(R(Z,Y)Z,X)-\Omega(Z,R(Z,Y)X)\\[8pt]
=g(R(Z,Y)Z,JX)+g(R(Z,X)Y,JZ).
\end{array}
$$
Hence
$$
S(\Omega)(X,Y)=Ricci(Y,JX)-Ricci(X,JY)+2\rho^{\ast}(X,JY).
$$
By (\ref{rho-star}), in order to show that $\rho^{\ast}(X,Y)=\rho^{\ast}(Y,X)$ for every $X,Y$, it
is enough to check that $\rho^{\ast}(X,Y)=\rho^{\ast}(Y,X)$ for all unit vectors $X,Y\in TM$ with
$g(X,Y)=g(X,JY)=0$. If $X$, $Y$ are such vectors,  $E_1=X$, $E_2=JE_1$, $E_3=Y$, $E_4=JY$ is an
orthonormal basis and the condition $g(S(\Omega),\tau)=0$ is equivalent to
\begin{equation}\label{S}
S(\Omega)(E_1,E_3)+S(\Omega)(E_4,E_2)=0, \quad S(\Omega)(E_1,E_4)+S(\Omega)(E_2,E_3)=0.
\end{equation}
These identities are equivalent to $\rho^{\ast}(E_1,E_4)=\rho^{\ast}(E_4,E_1)$ and
$\rho^{\ast}(E_1,E_3)=\rho^{\ast}(E_2,E_4)$ where
$\rho^{\ast}(E_2,E_4)=\rho^{\ast}(JE_1,JE_3)=\rho^{\ast}(E_3,E_1)$. Taking into account
(\ref{rho-star}) we see that (\ref{S}) is equivalent to $\rho^{\ast}(X,Y)=\rho^{\ast}(Y,X)$.

\section{Examples}

\subsection {Generalized Hopf surfaces}

Clearly, if $M$ is locally conformally K\"ahler ($d\theta=0$) and the Lee form $\theta$ has
constant length, the hypersurface $\Sigma_J$ is minimal. We have $||\theta||\equiv const$ on every
homogeneous locally conformally K\"ahler manifold. Also, if $\theta$ is parallel, then $d\theta=0$
and $||\theta||\equiv const$. Recall that a Hermitian surface with parallel Lee form is called a
generalized Hopf surface \cite{V} (or a Vaisman surface \cite{DO}); we refer to \cite{DO, V} for
basic properties and examples of such surfaces. The product of a Sasakian $3$-manifold with ${\Bbb
R}$ or $S^1$ admits a structure of generalized Hopf surface in a natural way. Conversely, every
such a surface locally is the product of a Sasakian manifold and  ${\Bbb R}$ \cite{V} (cf. also
\cite{GO}). A global structure theorem for compact generalized Hopf manifolds is obtained in
\cite{LV}.

As it is shown in \cite{Tr}, certain Inoue surfaces admit locally conformally K\"ahler structures
with $||\theta||\equiv const$ and non-parallel Lee form $\theta$.

\smallskip

Fix two complex numbers $\alpha$ and $\beta$ such that $|\alpha|\geq |\beta|>1$. Let
$\Gamma_{\alpha,\beta}$ be the group of transformation of ${\Bbb C}^2\setminus\{0\}$ generated by
the transformation $(u,v)\to (\alpha u,\beta v)$. Then, by a result of \cite{GO}, the quotient
$M_{\alpha,\beta}=({\Bbb C}^2\setminus\{0\})/\Gamma_{\alpha,\beta}$ admits a structure of a
generalized Hopf surface. Note that $M_{\alpha,\beta}$ (as any primary Hopf complex surface) is
diffeomorphic to $S^3\times S^1$.

\medskip

\noindent

{\bf The sypersurface $\Sigma_J$ in the twistor space of $S^3\times S^1$}

 We shall consider the Hopf surface $S^3\times S^1$ with its standard complex structure
${\cal J}$ and the product metric.

 According to \cite[Example 5]{D1}, we have
$$
{\cal Z}_+(S^3\times S^1)\cong\{[z_1,z_2,z_3,z_4]\in{\Bbb C}{\Bbb
P}^3:~|z_1|+|z_2|=|z_3|+|z_4|\}\times S^1.
$$
In order to give an explicit description of this isomorphism, we first recall that the twistor
space of an odd-dimensional oriented Riemannian manifold $(M,g)$ is the bundle ${\cal C}_+(M)$ over
$M$ whose fibre at a point $p\in M$ consists of all (linear) contact structures on the tangent
space $T_pM$ compatible with the metric and the  orientation, i.e. pairs $(\varphi,\xi)$ of
endomorphism $\varphi$ of $T_pM$ and a unit vector $\xi\in T_pM$ such that
$\varphi^2X=-X+g(X,\xi)\xi$,\, $g(\varphi X,\varphi Y)=g(X,Y)-g(X,\xi)g(Y,\xi)$ for $X,Y\in T_pM$,
and the orientation of $T_pM$ is induced by the orthogonal decomposition $T_pM=Im\,\varphi\oplus
{\Bbb R}\xi$ where the vector space $Im\,\varphi$ is oriented by means of the complex structure
$\varphi|Im\,\varphi$ on it. We refer to \cite{D1,D2} and the references therein for more
information about the twistor spaces of odd-dimensional manifolds. The twistor space ${\cal
C}_+(M)$ admits a $1$-parameter family of Riemannian metrics $h^c_t$ defined in a way similar to
the definition of the metrics $h_t$ on the twistor space ${\cal Z}_+$.

As is well-known, given $(\varphi,\xi)$ and $a\in S^1$, we can define a complex structure $I$ on
$T_pM\times T_aS^1$ in the following way.  Denote by $\displaystyle{\frac{\partial}{\partial t}}$
the vector field on $S^1$ determined by the local coordinate $e^{it}\to t$. Then set $I=\varphi$ on
$Im\,\varphi$, $I\xi=\displaystyle{(\frac{\partial}{\partial t})(a)}$,
$I\displaystyle{(\frac{\partial}{\partial t})(a)}=-\xi$. The complex structure $I$ is compatible
with the product metric of $M\times S^1$ and its orientation, $S^1$ (as well as any other sphere)
being oriented by the inward normal vector field. In this way we have a map
$$
F:{\cal C}_+(M)\times S^1\to {\cal Z}_+(M\times S^1).
$$
Endow ${\cal C}_+(M)\times S^1$ with the product metric. It is a simple observation that the map
$F$ is a bundle isomorphism preserving the metrics (and having other nice properties) \cite[Example
4]{D1}. Now we apply this observation to the case when $M=S^3$ and shall define an embedding of
${\cal C}_+(S^3)$ into ${\Bbb C}{\Bbb P}^3$ as in \cite[Examples 2 and 3]{D2}.

Denote the standard basis of ${\Bbb R}^6$ by $a_1,...,a_6$ and consider ${\Bbb R}^3$ and ${\Bbb
R}^4$ as the subspaces $span\{a_1,a_2,a_3\}$ and $span\{a_1,...,a_4\}$. Let $(\varphi,\xi)\in {\cal
C}_+(S^3)$ with $\varphi\in End(T_pS^3)$ and $\xi\in T_pM$, $p\in S^3$. Then we define a complex
structure $J$ on ${\Bbb R}^6$ by means of the orthogonal decomposition
$$
{\Bbb R}^6=Im\,\varphi\oplus{\Bbb R}\xi\oplus{\Bbb R}\{-p\}\oplus{\Bbb R}a_5\oplus{\Bbb R}a_6
$$
setting $J=\varphi$ on $Im\,\varphi$, $J\xi=-a_5$, $Jp=-a_6$, $Ja_5=\xi$, $Ja_6=p$. In this way we
obtain an embedding $\kappa$ of ${\cal C}_+(S^3)$ into the space $J_+({\Bbb}R^6)$ of complex
structures on ${\Bbb R}^6$ compatible with the metric and the orientation. The tangent space of
$J_+({\Bbb}R^6)$ at any point $I$ consists of skew-symmetric endomorphisms $Q$ of ${\Bbb R}^6$
anti-commuting with $I$. Denote by $G$ the standard metric $-\frac{1}{2}Trace\,PQ$ on the space of
skew-symmetric endomorphisms. Then $\frac{1}{2}\kappa^{\ast}G=h^c_{1/2}$ \cite[Example 2]{D2}. It
is well-known that $J_+({\Bbb}R^6)$ and ${\Bbb C}{\Bbb P}^3$ are both isomorphic to the twistor
space of $S^4$ (see, for example, \cite{Will}), so $J_+({\Bbb}R^6)\cong{\Bbb C}{\Bbb P}^3$; for a
direct proof see \cite{AGS,W}). We shall make use of the biholomorphism that sends a point
$[z_1,z_2,z_3,z_4]\in{\Bbb C}{\Bbb P}^3$ to the complex structure $J$ of ${\Bbb R}^6$ defined as
follows: Let $a_1,...,a_6$ be the standard basis of ${\Bbb R}^6$ and set $A_k=
\displaystyle{\frac{1}{\sqrt 2}}(a_{2k-1}-ia_{2k}), k=1,2,3$. Then the structure $J$ is given by
$$
-i|z|^2JA_1=(|z_1|^2-|z_2|^2-|z_3|^2+|z_4|^2)A_1+2\overline z_1 z_2A_2 + 2\overline z_1 z_3A_3+
2\overline z_4 z_3\overline A_2-2\overline z_4 z_2\overline A_3
$$
$$
 -i|z|^2JA_2=2\overline z_2 z_1
A_1+(-|z_1^2+|z_2|^2-|z_3|^2+|z_4|^2)A_2+ 2\overline z_2 z_3A_3- 2\overline z_4 z_3\overline
A_1+2\overline z_4 z_1\overline A_3
$$
$$
 -i|z|^2JA_3=2\overline z_3 z_1A_1+2\overline z_3
z_2A_2+(-|z_1|^2-|z_2|^2+|z_3|^2+|z_4|^2)A_3 +2\overline z_4 z_2\overline A_1-2\overline z_4
z_1\overline A_2.
$$
where $z=(z_1,z_2,z_3,z_4)$.

  For every $p\in S^3$, denote by $\times$ the vector cross-product on the oriented $3$-dimensional Euclidean space
$T_pS^3$. If $(\varphi,\xi)$ is a linear contact structure on $T_pM$ compatible with the metric and
the orientation, then
$$
\varphi (v)=\xi\times v, \quad v\in T_pS^3.
$$
In particular, $(\varphi,\xi)$ is uniquely determined by $\xi$. Define an oriented orthonormal
global frame of the bundle $TS^3$ by
$$
\xi_1(p)=(-p_2,p_1,-p_4,p_3),\quad \xi_2(p)=(-p_3,p_4,p_1,-p_2),\quad \xi_3(p)=(-p_4,-p_3,p_2,p_1),
$$
$$
\mbox {for}~(p_1,p_2,p_3,p_4)\in S^3.
$$
Set $\varphi_1(v)=\xi_1(p)\times v$, $v\in T_pS^3$. Then the standard complex structure ${\cal J}$
of $S^3\times S^1$ corresponds to the section $(\varphi_1,\xi_1)\times
\displaystyle{\frac{\partial}{\partial t}}$ of ${\cal C}_+(S^3)\times S^1$ under the isomorphism
$F: {\cal C}_+(S^3)\times S^1\to {\cal Z}_+(S^3\times S^1)$. We note also that if $J'$, $J''$
corresponds to $(\varphi',\xi')$, $(\varphi'',\xi'')$ under $F$, then $G(J',J'')=g(\xi',\xi'')$. In
particular, $J'$ and $J''$ are orthogonal if and only if $\xi'$ and $\xi''$ are so. Let
$(\varphi,\xi)\in {\cal C}_+(S^3)$ and $\xi\perp \xi_1(p)$, thus
$\xi=\lambda_2\xi_2(p)+\lambda_3\xi_3(p)$ where $\lambda_2^2+\lambda_3^2=1$. Then the point
$[z_1,...,z_4]\in{\Bbb C}{\Bbb P}^3$ corresponding to $(\varphi,\xi)$ under the embedding ${\cal
C}_+(S^3)\hookrightarrow J_+{\Bbb R}^6\cong {\Bbb C}{\Bbb P}^3$ is given by
$$
\begin{array}{c}
z_1=\displaystyle{\frac{1}{2}[-(p_1+ip_2)-(\lambda_3-i\lambda_2)(p_3-ip_4)],
z_2=\frac{1}{2}[(\lambda_3-i\lambda_2)(p_1-ip_2)-(p_3+ip_4)},\\[8pt]
\displaystyle{z_3=\frac{1}{2}, z_4=-\frac{1}{2}(\lambda_3-i\lambda_2)}.
\end{array}
$$
In particular, we have $4|z_3|^2=4|z_4|^2=|z|^2$. Conversely, let $[z_1,...,z_4]\in{\Bbb C}{\Bbb
P}^3$ be a point for which $4|z_3|^2=|z|^2$ and $4|z_4|^2=|z|^2$. Let $p_1,...,p_4$,
$\lambda_2,\lambda_3$ be the real numbers determined by the equations
$$
p_1+ip_2=-\frac{2(z_1\bar z_3+\bar z_2 z_4)}{|z|^2},\> p_3+ip_4=\frac{2(\bar z_1 z_4-z_2 \bar
z_3)}{|z|^2},\> \lambda_3+i\lambda_2=-\frac{4z_3\bar z_4}{|z|^2}.
$$
Then $p=(p_1,...,p_4)\in S^3$, $\lambda_2^2+\lambda_3^2=1$ and $[z_1,...,z_4]$ corresponds under
the embedding ${\cal C}_+(S^3)\hookrightarrow J_+{\Bbb R}^6\cong {\Bbb C}{\Bbb P}^3$ to
$(\varphi,\xi)$  determined by $\xi=\lambda_2\xi_2(p)+\lambda_3\xi_3(p)$. It follows that
$$
\Sigma_{{\cal J}}\cong \{[z_1,z_2,z_3,z_4]\in {\Bbb C}{\Bbb P}^3:~ 4|z_3|^2= 4|z_4|^2=|z|^2\}\times
S^1.
$$

\subsection {Kodaira surfaces}

Recall that every primary Kodaira surface $M$ can be obtained in the following way
\cite[p.787]{Kodaira}. Let $\varphi_k(z,w)$ be the affine transformations of ${\Bbb C}^2$ given by
$$\varphi_k(z,w) = (z+a_k,w+\overline{a}_kz+b_k),$$ where $a_k$, $b _k$, $k=1,2,3,4$, are complex
numbers such that
$$
a_1=a_2=0, \quad Im(a_3{\overline a}_4) =m b_1\neq 0,\quad b_2\neq 0
$$
for some integer $m>0$. They generate a group $G$ of transformations acting freely and properly
discontinuously on ${\Bbb C}^2$, and $M$ is the  quotient space ${\Bbb C}^2/G$.

It is well-known that $M$ can also be describe  as the quotient of ${\Bbb C}^2$ endowed with a
group structure by a discrete subgroup $\Gamma$. The multiplication on ${\Bbb C}^2$ is defined by
$$
(a,b).(z,w)=(z+a,w+\overline{a}z+b),\quad (a,b), (z,w)\in  {\Bbb C}^2,
$$
and $\Gamma$ is the subgroup generated by $(a_k,b_k)$, $k=1,...,4$ (see, for example, \cite{Borc}).

Further we shall consider $M$ as the quotient of the group ${\Bbb C}^2$ by the discrete subgroup
$\Gamma$. Every left-invariant object on ${\Bbb C}^2$ descends to a globally defined object on $M$
and both of them will be denoted by the same symbol.

We identify ${\Bbb C}^2$ with ${\Bbb R}^4$ by $(z=x+iy,w=u+iv)\to (x,y,u,v)$ and set
$$
A_1=\frac{\partial}{\partial x}-x\frac{\partial}{\partial u}+y\frac{\partial}{\partial v},\quad
A_2=\frac{\partial}{\partial y}-y\frac{\partial}{\partial u}-x\frac{\partial}{\partial v},\quad
A_3=\frac{\partial}{\partial u},\quad A_4=\frac{\partial}{\partial v}.
$$
These fields form a basis for the space of left-invariant vector fields on ${\Bbb C}^2$. We note
that the Lie brackets of the vector fields $A_1,...,A_4$ are
$$
[A_1,A_2]=-2A_4,\quad [A_i,A_j]=0\>\mbox{for all other}\> i,j,\> i<j.
$$
It follows that the group ${\Bbb C}^2$ defined above is solvable.

 Denote by $g$ the left-invariant Riemannian metric on $M$ for which the basis
$A_1,...,A_4$ is orthonormal.

We shall consider almost complex structures $J$ on $M$ compatible with the metric $g$  obtained
from left-invariant almost complex structures on ${\Bbb C}^2$. Note that by \cite{Has} every
complex structure on $M$ is induced by a left-invariant complex structure.

\noindent {\bf I}. If $J$ is a left-invariant almost complex structure compatible with $g$, we have
$JA_i=\sum_{j=1}^4 a_{ij}A_j$ where $a_{ij}$ are constants with $a_{ij}=-a_{ji}$.  Let $N$ be the
Nijenhuis tensor of $J$. Computing $N(A_i,A_j)$ in terms of $a_{ij}$, one can see \cite{M} that $J$
is integrable if and only if
$$
JA_1=\varepsilon_1 A_2,\quad JA_3=\varepsilon_2 A_4,\quad \varepsilon_1,\varepsilon_2=\pm 1.
$$
Since we are dealing  with the complex structures orthogonal to $J$, it is enough to consider the
two structures $J_{\varepsilon}$ defined by
$$
JA_1=\varepsilon A_2,\quad JA_3=A_4,\quad \varepsilon=\pm 1.
$$
Endow  $M$ with the orientation induced by $J_{\varepsilon}$. Then $\Lambda^2_{+}M$ admits a global
orthonormal frame defined by
$$
s_1^{\varepsilon}=\varepsilon A_1\wedge A_2+A_3\wedge A_4,\quad s_2^{\varepsilon}=A_1\wedge
A_3+\varepsilon A_4\wedge A_2,\quad s_3^{\varepsilon}=A_1\wedge A_4+\varepsilon A_2\wedge A_3.
$$
Hence we have a natural diffeomorphism
$$
F^{\varepsilon}: {\cal Z}_{+}(M)\cong M\times S^2, \quad \sum_{k=1}^3 x_k s_k^{\varepsilon}(p)\to
(p,x_1,x_2,x_3),
$$
under which
$$
\Sigma_{J_{\varepsilon}}\cong \{(p,x)\in M\times S^2:~x_1=0\}.
$$

In order to find an explicit formula for the metrics $h_t$ we compute the covariant derivatives of
$s_1^{\varepsilon},s_2^{\varepsilon},s_3^{\varepsilon}$ with respect to the Levi-Civita connection
$\nabla$ of $g$. The non-zero covariant derivatives $\nabla_{A_i}A_j$ are
$$
\nabla_{A_1}A_2=-\nabla_{A_2}A_1=-A_4,\quad \nabla_{A_1}A_4=\nabla_{A_4}A_1=A_2,\quad
\nabla_{A_2}A_4=\nabla_{A_4}A_2=-A_1.
$$
Then
$$
\nabla_{A_1}s_1^{\varepsilon}=-\nabla_{A_4}s_2^{\varepsilon}=-\varepsilon s_3^{\varepsilon},\quad
\varepsilon \nabla_{A_1}s_3^{\varepsilon}=-\nabla_{A_2}s_2^{\varepsilon}=s_1^{\varepsilon};\quad
\nabla_{A_2}s_1^{\varepsilon}=-\varepsilon\nabla_{A_4}s_3^{\varepsilon}=s_2^{\varepsilon}.
$$
and all other covariant derivatives $\nabla_{A_i}s_k^{\varepsilon}$ are zero. It follows that
$F_{\ast}^{\varepsilon}$ sends the horizontal lifts $A_1^h,...,A_4^h$ at a point
$\sigma=\sum_{k=1}^3 x_k s_k^{\varepsilon}(p)\in{\cal Z}_{+}(M)$ to the vectors
$$
A_1+\varepsilon (-x_3,0,x_1),\quad A_2+(x_2,-x_1,0),\quad A_3,\quad A_4+\varepsilon (0,x_3,-x_2).
$$
For $x=(x_1,x_2,x_3)\in S^2$, set
$$
u_1^{\varepsilon}(x)=\varepsilon (-x_3,0,x_1),\quad u_2^{\varepsilon}(x)=(x_2,-x_1,0),\quad
u_3^{\varepsilon}(x)=0,\quad u_4^{\varepsilon}(x)=\varepsilon (0,x_3,-x_2).
$$
These are tangent vectors to $S^2$ at the point $x$.  Denote by $h_t^{\varepsilon}$ the pushforward
of the metric $h_t$ by $F^{\varepsilon}$ Then, if $X,Y\in T_pM$ and $P,Q\in T_xS^2$,
\begin{equation}\label{ht}
h_t^{\varepsilon}(X+P,Y+Q)=g(X,Y)+t<P-\sum_{i=1}^4 g(X,A_i)u_i^{\varepsilon}(x),Q-\sum_{j=1}^4
g(Y,A_j)u_j^{\varepsilon}(x)>
\end{equation}
where $<.,.>$ is the standard metric of ${\Bbb R}^3$.

Let $\theta_{\varepsilon}$ be the Lee form of the Hermitian manifold $(M,g,J_{\varepsilon})$. We
have $\theta_{\varepsilon}(X)=-2\varepsilon g(X,A_3)$ which implies $\nabla\theta_{\varepsilon}=0$. Therefore,
by Proposition~\ref{min}, the hypersurface $\{(p,x)\in M\times S^2:~x_1=0\}$ in $M\times S^2$ is
minimal with respect to the metrics $h_t^{\varepsilon}$ given by (\ref{ht}).

\smallskip

\noindent {\bf II}. Suppose again that $J$ is an almost complex structure on $M$ obtained from a
left-invariant almost complex structure on $G$ and compatible with the metric $g$. Denote  the
fundamental $2$-form of the almost Hermitian structure $(g,J)$ by $\Omega$. Set $JA_i=\sum_{j=1}^4
a_{ij}A_j$. The basis dual to $A_1,...,A_4$ is $\alpha_1=dx$, $\alpha_2=dy$, $\alpha_3=xdx+ydy+du$,
$\alpha_4=-ydx+xdy+dv$. We have $d\alpha_1=d\alpha_2=d\alpha_3=0$, $d\alpha_4=2dx\wedge dy$. Hence
$d\Omega=d\sum_{i<j}a_{ij}\alpha_i\wedge\alpha_j=-2a_{34}dx\wedge dy\wedge du$. Thus $d\Omega=0$ is
equivalent to $a_{34}=0$. If $a_{34}=0$, we have $a_{1j}=a_{2j}=0$ for $j=1,2$, $a_{3k}=a_{4k}=0$
for $k=3,4$, $a_{13}^2+a_{14}^2=1$, $a_{13}a_{23}+a_{14}a_{24}=0$, $a_{23}^2+a_{24}^2=1$. It
follows that  the structure $(g,J)$ is almost K\"ahler (symplectic) if and only if $J$ is given by
(\cite{M})
$$
\begin{array}{c}
JA_1=-\varepsilon_1\sin\varphi A_3+\varepsilon_1\varepsilon_2\cos\varphi A_4,\quad
JA_2=-\cos\varphi A_3-\varepsilon_2\sin\varphi A_4, \\[6pt]
JA_3=\varepsilon_1\sin\varphi A_1+\cos\varphi A_2,\quad JA_4=-\varepsilon_1\varepsilon_2\cos\varphi
A_1+\varepsilon_2\sin\varphi A_2, \\[6pt]
\varepsilon_1,\varepsilon_2=\pm 1, \> \varphi\in{\Bbb R}.
\end{array}
$$
For fixed  $\epsilon=(\varepsilon_1,\varepsilon_2)$ and $\varphi$, denote by $J^{\epsilon,\varphi}$
the almost complex structure defined by these identities.  Set
$$
E_1=A_1,\quad E_2=-\varepsilon_1\sin\varphi A_3+\varepsilon_1\varepsilon_2\cos\varphi A_4,\quad
E_3=\cos\varphi A_3+\varepsilon_2\sin\varphi A_4,\quad E_4=A_2.
$$
Then $E_1,...,E_4$ is an orthonormal frame of $TM$ for which $J^{\epsilon,\varphi}E_1=E_2$ and
$J^{\epsilon,\varphi}E_3=E_4$. The only non-zero Lie bracket of these fields is
$$
[E_1,E_4]=-2(\varepsilon_1\varepsilon_2\cos\varphi E_2+\varepsilon_2\sin\varphi E_3).
$$
The non-zero covariant derivatives $\nabla_{E_i}E_j$ are
$$
\begin{array}{c}
\nabla_{E_1}E_2=\nabla_{E_2}E_1=\varepsilon_1\varepsilon_2\cos\varphi E_4,\quad
\nabla_{E_1}E_3=\nabla_{E_3}E_1=\varepsilon_2\sin\varphi E_4,\\[6pt]
\nabla_{E_1}E_4=-\nabla_{E_4}E_1=-\varepsilon_1\varepsilon_2\cos\varphi
E_2-\varepsilon_2\sin\varphi
E_3,\\[6pt]
\nabla_{E_2}E_4=\nabla_{E_4}E_2=-\varepsilon_1\varepsilon_2\cos\varphi E_1,\quad
\nabla_{E_3}E_4=\nabla_{E_4}E_3=-\varepsilon_2\sin\varphi E_1.
\end{array}
$$
Using  $E_1,...,E_4$, we define a global orthonormal frame $s_1^{+},s_2^{+},s_3^{+}$ of
$\Lambda^2_{+}TM$ via (\ref{s-basis}). We have the following table for the covariant derivatives of
$s_1^{+},s_2^{+},s_3^{+}$:
$$
\begin{array}{c}
\nabla_{E_1}s_1^{+}=\nabla_{E_4}s_2^{+}=\varepsilon_1\varepsilon_2\cos\varphi s_3^{+},\quad
\nabla_{E_4}s_1^{+}=-\nabla_{E_1}s_2^{+}=-\varepsilon_2\sin\varphi s_3^{+},\\[6pt]
\nabla_{E_2}s_1^{+}=\varepsilon_1\varepsilon_2\cos\varphi s_2^{+},\quad
\nabla_{E_2}s_2^{+}=-\varepsilon_1\varepsilon_2\cos\varphi s_1^{+},\quad
\nabla_{E_2}s_3^{+}=0,\\[6pt]
\nabla_{E_3}s_1^{+}=\varepsilon_2\sin\varphi s_2^{+},\quad
\nabla_{E_3}s_2^{+}=-\varepsilon_2\sin\varphi s_1^{+},\quad \nabla_{E_3}s_3^{+}=0,\\[6pt]
\nabla_{E_1}s_3^{+}=-\varepsilon_1\varepsilon_2\cos\varphi s_1^{+}-\varepsilon_2\sin\varphi
s_2^{+},\quad \nabla_{E_4}s_3^{+}=\varepsilon_2\sin\varphi s_1^{+}
-\varepsilon_1\varepsilon_2\cos\varphi s_2^{+}.
\end{array}
$$
The frame $s_1^{+},s_2^{+},s_3^{+}$ gives rise to an obvious diffeomeorphism $F^{\epsilon,\varphi}:
{\cal Z}_{+}(M)\cong M\times S^2$ for which $\Sigma_{J^{\epsilon,\varphi}}\cong \{(p,x)\in M\times
S^2:~x_1=0\}$. To describe the pushforward $h_t^{\varepsilon,\varphi}$ of the metric $h_t$ by
$F^{\varepsilon,\varphi}$, we set
$$
\begin{array}{c}
u_1^{\epsilon,\varphi}(x)=(x_3\varepsilon_1\varepsilon_2\cos\varphi, x_3\varepsilon_2\sin\varphi,-x_1\varepsilon_1\varepsilon_2\cos\varphi-x_2\varepsilon_2\sin\varphi),\\[6pt]
u_2^{\epsilon,\varphi}(x)=(x_2\varepsilon_1\varepsilon_2\cos\varphi,
-x_1\varepsilon_1\varepsilon_2\cos\varphi, 0) ,\quad
u_3^{\epsilon,\varphi}(x)=(x_2\varepsilon_2\sin\varphi,-x_1\varepsilon_2\sin\varphi,0)\\[6pt]
u_4^{\epsilon,\varphi}(x)=(-x_3\varepsilon_2\sin\varphi,
x_3\varepsilon_1\varepsilon_2\cos\varphi,x_1\varepsilon_2\sin\varphi-x_2\varepsilon_1\varepsilon_2\cos\varphi).
\end{array}
$$
for $x=(x_1,x_2,x_3)\in S^2$. Then, if $X,Y\in T_pM$ and $P,Q\in T_xS^2$,
\begin{equation}\label{ht-sympl}
h_t^{\epsilon,\varphi}(X+P,Y+Q)=g(X,Y)+t<P-\sum_{i=1}^4
g(X,E_i)u_i^{\epsilon,\varphi}(x),Q-\sum_{j=1}^4 g(Y,E_j)u_j^{\epsilon,\varphi}(x)>.
\end{equation}
It is easy to compute that
$$
\rho^{\ast}(E_1,E_4)=\rho^{\ast}(E_4,E_1)=-\varepsilon_1\sin\varphi.\cos\varphi,\quad
\rho^{\ast}(E_1,E_3)=\rho^{\ast}(E_3,E_1)=0.
$$
It follows from Proposition~\ref{min-sympl} that $\{(p,x)\in M\times S^2:~x_1=0\}$ is a minimal
hypersurface in  $M\times S^2$, the latter manifold being endowed with the metrics
$h_t^{\epsilon,\varphi}$ given by (\ref{ht-sympl}).

Secondary Kodaira surfaces are quotients of primary ones by groups of order 2,3,4 or 6. Every
secondary Kodaira surface is a homogeneous manifold (in fact a solvmanifold). It admits a basis of
left-invariant vector fields $A_1,...,A_4$  such that the complex structure sends $A_1$, $A_3$ to
$A_2$,$A_4$ and $[A_1,A_2]=-2A_4$, $2[A_3,A_1]=A_2$, $2[A_3,A_2]=-A_1$  \cite{Has}. Computations as
above show that if $J$ is a left-invariant complex or symplectic structure, then $\Sigma_J$ is a
minimal hypersurface.

\end{document}